\documentclass[article]{elsart}

\usepackage{amsmath}
\usepackage{amssymb}
\usepackage{mathrsfs}

     \DeclareMathOperator{\Aut}{Aut}

    \DeclareMathOperator{\Mod}{Mod}

    \DeclareMathOperator{\Char}{Char}
    \DeclareMathOperator{\Ext}{Ext} \DeclareMathOperator{\II}{II}
    \DeclareMathOperator{\vN}{vN}\DeclareMathOperator{\SL}{SL}
    \DeclareMathOperator{\III}{III}\DeclareMathOperator{\HT}{HT}
    
    \DeclareMathOperator{\gp}{\bf GP}
    \DeclareMathOperator{\gf}{\bf GF}
    \DeclareMathOperator{\wicc}{{\bf wT}_{ICC}}
    \DeclareMathOperator{\Out}{Out}
    \DeclareMathOperator{\itpf1}{ITPFI}
    \DeclareMathOperator{\Ad}{Ad}

    \newcommand{\act}{\mathcal A}
    \newcommand{\actson}{\curvearrowright}

\def\R{{\mathbb R}}
\def\C{{\mathbb C}}

\def\N{{\mathbb N}}
\def\Z{{\mathbb Z}}

\def\F{{\mathbb F}}

\def\T{{\mathbb T}}

\usepackage{mathrsfs}

\begin{document}

\begin{frontmatter}



\title{The classification problem for von Neumann factors}


\author[rs]{Roman Sasyk}
\author[at]{Asger T\"ornquist\corauthref{cor1}}

\address[rs]{Instituto de Ciencias and Instituto Argentino de Matem\'aticas-CONICET, Buenos Aires, Argentina}

\address[at]{Kurt Goedel Research Center, Vienna, Austria}

\corauth[cor1]{W\"ahringer Strasse 25, 1090 Vienna, Austria}

\ead{asger@logic.univie.ac.at}

\begin{abstract}
We prove that it is not possible to classify separable von Neumann
factors of types $\II_1$, $\II_\infty$ or $\III_\lambda$, $0\leq
\lambda\leq1$, up to isomorphism by a Borel measurable assignment of
``countable structures'' as invariants. In particular the
isomorphism relation of type $\II_1$ factors is not smooth. We also
prove that the isomorphism relation for von Neumann $\II_1$ factors
is analytic, but is not Borel.

\end{abstract}


\end{frontmatter}

\section{Introduction}

The purpose of this paper is to apply the notion of Borel
reducibility of equivalence relations, developed extensively in
descriptive set theory in recent years, and the deformation-rigidity
techniques of Sorin Popa, to study the global structure of the set
of factors on a separable Hilbert space.

Recall that if $E, F$ are equivalence relations on standard Borel
spaces $X$ and $Y$, respectively, we say that $E$ is {\it Borel
reducible} to $F$, written $E\leq_B F$, if there is a Borel $f:X\to
Y$ such that
$$
xEy\iff f(x)Ff(y),
$$
in other words, if it is possible to classify the points of $X$ up
to $E$ equivalence by a Borel assignment of invariants that are
$F$-equivalence classes. We write $E<_B F$ if $E\leq_B F$ but not
$F\leq_B E$.

In this paper, the space of study will be the space $\vN(\mathcal
H)$ of von Neumann algebras on a separable Hilbert space $\mathcal
H$. Effros introduced in \cite{effros1},\cite{effros2}, a standard
Borel structure on this space and showed that the set of factors
$\mathscr F(\mathcal H)$ is a Borel subset of $\vN(\mathcal H)$, in
particular is a standard Borel space. Let $\simeq^{\vN(\mathcal H)}$
denote the isomorphism relation in $\vN(\mathcal H)$ and
$\simeq^{\mathscr F(\mathcal H)}$ the isomorphism relation in
$\mathscr F(\mathcal H)$. It is natural to ask if it is possible to
classify isomorphism classes of factors by an assignment of
invariants which are countable groups, countable graphs, countable
linear orders, countable fields, or other kinds of ``countable
structures'' type invariants. Using the Axiom of Choice the answer
is in principle `yes', given that there are at most continuum many
isomorphism classes in $\mathscr F$ and continuum many
non-isomorphic countable groups. However,
 this kind of ``classification'' is of no interest since the Axiom of
Choice provides us with no concrete way of computing the invariants.
Therefore, the correct question to ask is if $\simeq^{\mathscr F}$
is Borel reducible to the isomorphism relation of countable graphs,
groups, fields, etc.

To be more specific, one can naturally regard the Polish space
\begin{align*}
\gp=\{(f,e)\in \N^{\N\times\N}\times\N: (\forall i,j,k\in N)
(f(f(i,j),k)=f(i,f(j,k)))\\
\wedge (\forall i\in \N) f(i,e)=f(e,i)=i \wedge (\forall i)(\exists
l) f(i,l)=e\}
\end{align*}
as the set of all countable groups up to isomorphism. The set of
countable graphs may naturally be identified with the Polish space
$$
\gf=\{g\in\{0,1\}^{\N\times\N}: (\forall i,j) g(i,j)=g(j,i)\wedge
(\forall i) g(i,i)=0\}.
$$
Denote by $\simeq^{\gp}$ the isomorphism relation in $\gp$, and
$\simeq^{\gf}$ the isomorphism relation in $\gf$. Then we may ask if
$\simeq^{\mathscr F(\mathcal H)}$ is Borel reducible to
$\simeq^{\gp}$ or $\simeq^{\gf}$.

The appropriate general context to phrase the question is to ask if
there is a countable language $\mathcal L$, in the sense of model
theory, such that $\simeq^{\mathscr F(\mathcal H)}$ is Borel
reducible to the isomorphism relation $\simeq^{\Mod(\mathcal L)}$ in
the space $\Mod(\mathcal L)$ of countable models of $\mathcal L$
(see \cite{hjorth}). If this were the case then we would say that
$\simeq^{\mathscr F(\mathcal H)}$ is {\it classifiable by countable
structures}. This generalization encompasses all natural countable
structures type invariants, including the two classes $\gp$ and
$\gf$ above. Our first result is that the answer is no:

{\thm The isomorphism relation for factors $\simeq^{\mathscr
F(\mathcal H)}$ is not classifiable by countable structures. In
fact, the restriction of $\simeq^{\mathscr F(\mathcal H)}$ to
\emph{any} of the classes $\II_1$, $\II_\infty$ and $\III_\lambda,
0\leq\lambda\leq 1$, is not classifiable by countable structures. In
particular, the isomorphism relation of type $\II_1$ factors is not
smooth.}

\bigskip
Along these lines, in \S 4 we also prove:

{\thm The isomorphism relation for injective type $\III_0$ factors
is not classifiable by countable structures.}
\medskip

It should be noted that Woods proved in \cite{woods} that
isomorphism of ITPFI factors is not smooth in the classical sense of
Mackey, Effros and Glimm, by showing that $E_0$, defined on
$\{0,1\}^\N$ by
$$
xE_0 y\iff (\exists N)(\forall n\geq N) x(n)=y(n),
$$
is Borel reducible to $\simeq^{\itpf1}$. Our conclusion is much
stronger than this. Namely, by a classical result of Baer, $E_0$ is
Borel bi-reducible to the isomorphism relation for countable rank 1
torsion free Abelian groups, and so $E_0$ classes are essentially
countable structures type invariants. It is known that while
$E_0\leq_B\simeq^{\gp}$ and $E_0\leq_B\simeq^{\gf}$, it also holds
that $\simeq^{\gp}\not\leq_B E_0$ and $\simeq^{\gf}\not\leq_B E_0$.
The reader is referred to \cite{kechris3} for an introduction to the
subject of Borel reducibility.

\medskip

For our next result, denote by $\wicc$ the subset of $\gp$
consisting of countably infinite ICC groups with relative property
(T) over an infinite normal subgroup. In \S 5 we show that the
isomorphism relation $\simeq^{\wicc}$ is Borel complete for
countable structures, that is, for every countable language
$\mathcal L$ and every invariant Borel class $\mathscr C\subseteq
\Mod(\mathcal L)$ we have $\simeq^{\mathscr C}\leq_B
\simeq^{\wicc}$. Combining this with a deep result of Popa in
\cite{popa3} we obtain:

{\thm It holds that $\simeq^{\wicc}<_B\simeq^{\mathscr F_{\II_1}}$.
Thus $\simeq^{\mathscr F_{\II_1}}$ is Borel complete for countable
structures, but not classifiable by countable structures.}

\bigskip

It follows from Theorem 3 that, in particular,
$\simeq^{\gf}<_B\simeq^{\mathscr F_{\II_1}}$ and
$\simeq^{\gp}<_B\simeq^{\mathscr F_{\II_1}}$. Theorem 3 also has the
following consequence:

{\cor The isomorphism relation $\simeq^{\mathscr F}$, regarded as a
subset of $\mathscr F\times\mathscr F$, is complete analytic. In
particular, it is not a Borel set.}

\bigskip

{\it Acknowledgements.} Research for this paper was initiated during
the Fields Institute Thematic Program on Operator Algebras during
the Fall of 2007. We would like to thank the Fields Institute for
kind hospitality. We would also like to thank George Elliott, Ilijas
Farah and Thierry Giordano for helpful discussions. A special thanks
is due to Stefaan Vaes for suggesting some of the arguments in
section 3. A. T\"ornquist was supported in part by the Danish
Natural Science Research Council post-doctoral grant no.
272-06-0211.

\section{The finite case}

\paragraph{The space $\act(G,X,\mu)$.} Let $(X,\mu)$ be a standard Borel
probability space. Recall that $\Aut(X,\mu)$ is the group of measure
preserving transformations of $X$, which is a Polish group when
given the topology it inherits when it is naturally identified with
a closed subgroup of the unitary group of $L^2(X,\mu)$. If $G$ is a
countable group, we let
$$
\act(G,X,\mu)=\{\sigma\in\Aut(X,\mu)^G: (\forall g_1, g_2\in
G)\sigma(g_1g_2)=\sigma(g_1)\sigma(g_2)\}.
$$

It is easy to see that this is a closed subset of $\Aut(X,\mu)^G$
when the latter is given the product topology. We identify
$\act(G,X,\mu)$ with the space of measure preserving $G$-actions on
$X$.

\paragraph{The Effros Borel space.} Let $\mathcal H$ be a separable
Hilbert space. We let $\vN(\mathcal H)$ denote the weakly closed,
$^*$-closed subalgebras of $\mathcal B(\mathcal H)$, that is,
$\vN(\mathcal H)$ is the space of von Neumann algebras. Effros has
shown that this is a standard Borel space if equipped with the Borel
structure generated by the sets
$$
\{M\in \vN(\mathcal H): M\cap U\neq \emptyset \}
$$
where $U$ is a weakly open subset of $\mathcal B(\mathcal H)$, see
\cite{effros1}. There is a natural Polish topology on $\vN(\mathcal
H)$ called the {\it Effros-Mar\'echal} topology. The reader is
referred to the detailed study by Haagerup and Winsl\o w in
\cite{haagwin1}, \cite{haagwin2}.

\paragraph{The group-measure space construction.} Let $G$ be a
discrete group and $(X,\mu)$ a standard Borel probability space. Let
$\sigma\in\act(G,X,\mu)$ be a m.p. $G$-action on $X$. We define on
$L^2(G\times X)$ the unitary operators
$$
U_g^{\sigma}(f)(h,x)=f(g^{-1}h,\sigma(g^{-1})(x)).
$$
and for each $\varphi\in L^\infty(X)$ the multiplication operators
$$
L_{\varphi}(f)(h,x)=\varphi(x)f(h,x).
$$
We denote by $L^\infty(X)\rtimes_{\sigma} G$ the finite von Neumann
algebra in $\mathcal B(L^2(G\times X))$ generated by the $U_g$,
$g\in G$, and $L_\varphi$, $\varphi\in L^\infty(X)$.

It is well-known that if $\sigma$ is ergodic then
$L^\infty(X)\rtimes_{\sigma} G$ is a factor, moreover, if the action
is free then $L^\infty(X)$ is a {\it Cartan subalgebra} of
$L^\infty(X)\rtimes_{\sigma} G$. (Recall that $A$ is a Cartan
subalgebra of $M$ if $A$ is a MASA and $\mathcal N_M(A)''=\{u\in
\mathcal U(M):\, uAu^*=A\}''=M$). By results of Feldman and Moore
\cite{FM}, two free ergodic m.p. actions of possible different
groups $\sigma\in \act(G,X,\mu)$ and $\theta\in \act(H,Y,\nu)$ are
orbit equivalent if and only if their corresponding inclusions of
Cartan subalgebras $L^\infty(X)\subset L^\infty(X)\rtimes_{\sigma}
G$\,\,, $L^\infty(Y)\subset L^\infty(Y)\rtimes_{\theta} H$ are
isomorphic. Thus the study of orbit equivalence of m.p. group
actions can be translated into a problem regarding inclusions of
finite von Neumann algebras. It is worth mentioning that isomorphism
between group-measure space von Neumann algebras does not imply
isomorphism of the corresponding Cartan subalgebra inclusions (see
\cite{connesjones}). However, Popa's deformation-rigidity machinery
\cite{popa} assures that this is the case for the particular kinds
of groups and group actions we study in this paper.

{\lem The map $\sigma\mapsto L^\infty(X)\rtimes_\sigma G$ is a Borel
function from $\act(G,X,\mu)$ to $\vN(L^2(G\times X))$ when the
latter is given the Effros Borel structure.\label{basiclemma}}

{\pf By the Corollary to Theorem 2 in \cite{effros1}, it is enough
to find a countable family $f_n:\act(G,X,\mu)\to \mathcal B
(L^2(G\times X))$ of Borel functions such that for each
$\sigma\in\act(G,X,\mu)$ we have that $(f_n(\sigma))_{n\in\N}$ is
dense in $L^\infty(X)\rtimes_\sigma G$.

Let $D_n$ be an enumeration of the dyadic intervals, and let
$\chi_{D_n}$ be the characteristic functions. $L^\infty(X)$ is
separable in the weak topology and is generated by the functions
$\chi_{D_n}$. Thus $L^\infty(X)\rtimes_\sigma G$ is generated by
$L_{\chi_{D_n}}$ and $U^\sigma_g$, $g\in G$. For each $g\in G$ the
function
$$
\act(G,X\mu)\to \mathcal B(L^2(G\times X)):\sigma\mapsto U_g^\sigma
$$
is clearly Borel (in fact continuous) when $\mathcal B(L^2(G\times
X))$ is given the weak topology, and so are the constant functions
$$
\act(G,X,\mu)\to B(L^2(G\times X)):\sigma\mapsto L_{\chi_{D_n}}.
$$
Suppose now that $f,g:\act(G,X,\mu)\to \mathcal B(L^2(G\times X))$
are Borel functions. Then for all scalars $r_1,r_2\in\C$ we have
that the functions
$$
r_1f+r_2g, f^*, f\circ g
$$
are Borel. Hence the smallest class of functions containing
$\sigma\mapsto U^\sigma_g$ and $\sigma\mapsto L_{\chi_{D_n}}$ and
closed under taking linear combinations with rational-complex
scalars, adjoint and composition is a countable class of Borel
functions $f_n:\act(G,X,\mu)\to \mathcal B(L^2(G\times X))$,
$n\in\N$, and for each $\sigma\in\act(G,X,\mu)$ clearly
$(f_n(\sigma))_{n\in\N}$ is dense in $L^\infty(X)\rtimes_\sigma
G$.\qed }

Let $\sigma:\F_2\actson \T^2=X$ be the usual action on $\F_2$ on
$\T^2$, when $\F_2$ is viewed as a (finite index) subgroup of
$\SL(2,\Z)$. It was shown in \cite[\S 3]{at1} that there is a dense
$G_\delta$ set $\Ext(\sigma)\subseteq\Aut(X,\mu)$ extending the
action $\sigma$ to an a.e. free $\F_3$-action, specifically, if
$\F_2=\langle a,b\rangle$ and $T_a,T_b\in \Aut(X,\mu)$ are the
transformation corresponding to the generators, then the set
$$
\Ext(\sigma)=\{S\in\Aut(X,\mu): T_a,T_b \text{ and } S \text{
generate an a.e. free action of } \F_3\}
$$
is a dense $G_\delta$ set. For each $S\in\Ext(\sigma)$ we denote by
$\sigma_S$ the corresponding $\F_3$ action. It is easy to verify
that $S\mapsto\sigma_S$ is a continuous map $\Ext(\sigma)\to\mathcal
A(\F_3,X,\mu)$. It was shown in \cite[\S 5]{at1} that the relation
$$
S_1\sim_{oe} S_2\iff \sigma_{S_1} \text{ is orbit equivalent to }
\sigma_{S_2}
$$
has meagre classes, that there is a co-meagre set of dense classes,
and that $E_0\leq_B\sim_{oe}$. The following strengthening was also
noted by Kechris in \cite[Theorem 17.1]{kechris2}:

{\lem $\sim_{oe}$ is not classifiable by countable structures.}

{\pf Let $[E_{\sigma}]$ denote the full group of the $\F_2$-action
$\sigma$. As noted in \cite[p. 280]{at1}, the conjugation action of
$[E_{\sigma}]$ on $\Aut(X,\mu)$ preserves $\sim_{oe}$. Since
$\sigma$ is ergodic there is an ergodic m.p. $\Z$-action on $X$ such
that $E_\Z\subseteq E_\sigma$. By \cite[Theorem 11 and Claim
13]{foremanweiss}, the conjugation action of $[E_\Z]$ on
$\Aut(X,\mu)$ is turbulent. Hence $\sim_{oe}$ has meagre classes,
has a dense class and contains a turbulent action, so it follows
from \cite[Theorem 3.18]{hjorth} that it is generically
$S_\infty$-ergodic, in particular, it is not classifiable by
countable structures. \qed}

{\thm The von Neumann $\II_1$ factors that arise in the group
measure construction from an ergodic a.e. free $\F_3$ action are not
classifiable by countable structures.\label{basicthm}}

{\pf It suffices to show that for $S_1,\,S_2\in\Ext(\sigma)$,
$L^\infty(X)\rtimes_{\sigma_{S_1}}\F_3$ is isomorphic to
$L^\infty(X)\rtimes_{\sigma_{S_2}}\F_3$  if and only if
$\sigma_{S_1}$ and $\sigma_{S_2}$ are orbit equivalent. Namely, if
we can do this then by the previous Lemmas the map
$$
S\mapsto L^\infty(X)\rtimes_{\sigma_S} \F_3
$$
is a Borel reduction of $\sim_{oe}$ to isomorphism in $\mathscr
F_{\II_1}(L^2(\F_3\times X))$.

Since orbit equivalent actions give rise to isomorphic factors, we
only need to verify that the converse holds for the actions of the
form $\sigma_S$. By Feldman and Moore's Theorem cited above, it is
enough to show that an isomorphism of factors of the form
$L^\infty(X)\rtimes_{\sigma_S} \F_3$, $S\in\Ext(\sigma)$ gives an
isomorphism of the corresponding Cartan inclusions
$L^\infty(X)\subset L^\infty(X)\rtimes_{\sigma_S} \F_3$. This
follows from the seminal work of Popa \cite{popa} on $\II_1$ factors
with trivial fundamental group\footnote{We refer to the reader to
\cite{popa} for the pertinent definitions.}. Indeed, since $\F_3$
has the Haagerup compact approximation property, by proposition 3.1
in \cite{popa}, we have that for every $S\in \Ext(\sigma)$,
$L^\infty(X)\subset L^\infty(X)\rtimes_{\sigma_S} \F_3$ has the
relative property H. Since $L^\infty(X)\subset
L^\infty(X)\rtimes_\sigma\F_2$ is a rigid embedding and
$L^\infty(X)\rtimes_\sigma\F_2\subset L^\infty(X)\rtimes_{\sigma_S}
\F_3$, then by proposition 4.6 in \cite{popa}, $L^\infty(X)\subset
L^\infty(X)\rtimes_{\sigma_S}\F_3$ is a rigid embedding. Thus
$L^\infty(X)$ is an $HT_s$ Cartan subalgebra of
$L^\infty(X)\rtimes_{\sigma_S} \F_3$. But then by Corollary 6.5 in
\cite{popa}, any isomorphism between
$L^\infty(X)\rtimes_{\sigma_{S_1}} \F_3$ and
$L^\infty(X)\rtimes_{\sigma_{S_2}}\F_3$ can be perturbed by a
unitary to give an isomorphism of the HT Cartan subalgebras
inclusions. It follows that $\sigma_{S_1}$ and $\sigma_{S_2}$ are
orbit equivalent. \qed}

{\cor Let $\mathcal L$ be a countable language and let
$\Mod(\mathcal L)$ be the Polish space of countable models of
$\mathcal L$. Then it is not possible in Zermelo-Fraenkel set theory
\emph{without} the Axiom of Choice to construct a function
$$
f:\mathscr F_{\II_1}\to\Mod(\mathcal L)
$$
such that
$$
M_1\simeq^{\mathscr F_{\II_1}} M_2\iff f(M_1)\simeq^{\Mod(\mathcal
L)} f(M_2).
$$
In particular, it is not possible to construct such function with
codomain $\gp$ or $\gf$.}

{\pf By well-known Theorems of Solovay \cite{solovay} and Shelah
\cite{shelah}, it is consistent with ZF that every function from one
Polish space to another is Baire measurable. Hence if a function as
above could be constructed in ZF, then by the previous Theorem it
would be consistent that a Baire measurable function $\tilde
f:\Aut(X,\mu)\to\Mod(\mathcal L)$ existed such that
$$
S_1\sim_{oe} S_2\iff \tilde f(S_1)\simeq^{\Mod(\mathcal L)}\tilde
f(S_2).
$$
However, this is not possible since $\sim_{oe}$ is generically
$S_\infty$-ergodic. \qed}

{\it Remark.} Solovay and Shelah's Theorems are true even if we
allow the {\it countable} Axiom of Choice (countable AC) or the
principle of Dependent Choices (DC). Thus it holds that no function
as in the Corollary can be constructed in ZF+countable AC or ZF+DC.

\section{The semifinite and the purely infinite case}

In this section we prove the following:

{\thm The isomorphism relation for factors of types $\II_{\infty}$
and $\III_\lambda$, $0\leq\lambda\leq 1$,
are not classifiable by countable structures.}\\

Mimicking the argument given in the type $\II_1$ case, we exhibit
concrete examples of families of factors of types $\II_\infty$ and
$\III_\lambda$, $0\leq \lambda\leq 1$, for which $\sim_{oe}$ is
Borel  reducible to the isomorphism relation in each of the
families.  For the proof we rely on the main result of the previous
section, on a recent theorem of Popa concerning unique tensor
product decomposition of McDuff $\II_1$ factors and on Connes and
Takesaki cross product characterization of type $\III$ factors.
Throughout this section, for each  $S\in \Ext(\sigma)$ we denote by
$M_S$ the type $\II_1$ factor
$L^\infty(X,\mu)\rtimes_{\sigma_S}\F_3$ constructed in the previous
section, with $R$ the unique injective factor of type $\II_1$ and
with $R_{0,1}=R\otimes\mathcal B(\ell^2(\N))$ the unique type
$\II_\infty$ injective factor. We need this simple lemma, a proof of
which can be found in \cite[Corollary 3.8]{haagwin1}.

{\lem For each $N\in \vN(\mathcal H)$, the map
$\otimes_N:\vN(\mathcal H)\to\vN(\mathcal H\otimes \mathcal H)$
$M\to M\otimes N$ is Borel.}

\begin{pf}[Proof of Theorem 9:]

\medskip
\noindent {\it The type $\II_\infty$ case:} For each  $S\in
\Ext(\sigma)$,  $M_S\otimes\mathcal B(\ell^2(\N))$ is a factor of
type $\II_\infty$. Assume that $ \Theta$ is an isomorphism from
$M_{S_1}\otimes\mathcal B(\ell^2(\N))$ to $M_{S_2}\otimes\mathcal
B(\ell^2(\N))$. Then $M_{S_1}$ is isomorphic to an amplification of
$M_{S_2}$. Indeed, if $p={\rm Id}\otimes e_{11}$ then
$p\,M_{S_1}\otimes\mathcal B(\ell^2(\N))\,p\cong M_{S_1}$. Since
$\Theta(p)$ is a finite projection on the factor
$M_{S_2}\otimes\mathcal B(\ell^2(\N))$, there exist two finite
projections $c_1$ and $d_1$ in $M_{S_2}$ and $\mathcal
B(\ell^2(\N))$ such that $\Theta(p)\sim c_1\otimes d_1$. Since
$M_{S_2}\otimes\mathcal B(\ell^2(\N))$ is a properly infinite
factor, there exist two families of mutually equivalent and pairwise
orthogonal projections $(p_n)$ and $(q_n)$ such that $\sum p_n=\sum
q_n=1$, $p_1=\Theta(p)$, $q_1= c_1\otimes d_1$ and $p_n\sim q_n$.
Denote by $u_n\in M_{S_1}\otimes\mathcal B(\ell^2(\N))$ the partial
isometry such that $p_n=u_n^*u_n$ and $q_n=u_nu_n^*$. Then $u=\sum
u_n$ is a unitary and $u p_1 u^*=q_1$. Thus $\Ad(u)\circ \Theta$
implements an isomorphism between $M_{S_1}\cong
p\,M_{S_1}\otimes\mathcal B(\ell^2(\N))\,p $ and
$q_1\,M_{S_2}\otimes \mathcal B(\ell^2(\N))\,q_1=M_{S_2}^t$ where
$t=\tau_2\otimes {\rm Tr }\,(q_1)$. By Corollary 8.2 in \cite{popa}
the $\ell^2_{\HT}$ Betti numbers
 of the factors $M_S$, $S\in \Ext(\sigma)$, are all equal to
the Atiyah $\ell^2$-Betti number of $\F_3$, thus
$\beta^1_{\HT}(M_S)=\beta_1(\F_3)=2\neq 0$. Moreover,
$\beta^1_{\HT}(M_{S_2}^t)=\beta^1_{\HT}(M_{S_2})/t=2/t$. We conclude
that $t=1$. Thus the map $\Ext(\sigma)\ni S\mapsto
M_S\otimes\mathcal B(\ell^2(\N))$, is a Borel reduction from
$\sim_{oe}$ to isomorphism of $\II_\infty$ factors.

\bigskip
\noindent {\it The type $\III_\lambda$, $0<\lambda<1$, case:} Let
$R_\lambda$ , $0<\lambda< 1$, be the unique injective factor of type
$\III_\lambda$ and let $R_\lambda=R_{0,1}\rtimes_\theta\Z$ be its
discrete decomposition, where $R_{0,1}$ is the injective
$\II_\infty$ factor.
 We consider the family of type
$\III_\lambda$ factors $M_{S,\lambda}=M_S\otimes R_\lambda$, $S\in
\Ext(\sigma)$. Since $M_{S,\lambda}=M_S\otimes(
R_{0,1}\rtimes_\theta\Z)= (M_S\otimes
R_{0,1})\rtimes_{Id\otimes\theta}\Z$ , by the uniqueness of the
discrete decomposition \cite[Theorem XII.2.1]{takesaki2} two factors
 $M_{S_1,\lambda}$ and $M_{S_2,\lambda}$ are
isomorphic if and only if their cores $M_{S_1}\otimes R_{0,1}$ and
$M_{S_2}\otimes R_{0,1}$ are isomorphic. The same argument as in the
semifinite case shows that this occurs if and only if
$M_{S_1}\otimes R$ is isomorphic to an amplification of
$M_{S_2}\otimes R$, which in turn is isomorphic to $M_{S_2}\otimes
R$, since the fundamental group of $R$ is $\mathbb R_{+}^{*}$.

For every $S\in \Ext(\sigma)$, $M_S$ is a non-$\Gamma$ factor. By
Theorem 2 in \cite{popa2}, $M_{S_1}\otimes R$ is isomorphic to
$M_{S_2}\otimes R$ if and only if $M_{S_1}$ is isomorphic to an
amplification of $M_{S_2}$. As before,
 this occurs if and only if $M_{S_1}$ is isomorphic to $M_{S_2}$.
Thus the maps $\Ext(\sigma)\ni S\mapsto M_{S,\lambda}$,
$0<\lambda<1$, are
 Borel reductions from $\sim_{oe}$ to isomorphism of $\III_\lambda$ factors.

\bigskip
\noindent {\it The type $\III_1$ case:} The same argument used in
the type $\III_\lambda$ case, $0< \lambda < 1$, carries on to the
$\III_1$ case once we replace Connes discrete decomposition with
Takesaki continuous decomposition. Indeed, let  $R_\infty$ be the
unique injective factor of type $\III_1$. If $\phi$ is a faithful
semifinite normal weight on $R_\infty$, $\{\sigma^\phi\}$ denotes
the modular automorphism group associated to the weight $\phi$. By
Takesaki duality, $R_\infty= R_{0,1}\rtimes_{\sigma^\phi}\R$, where
$R_{0,1}$ is the injective $\II_\infty$ factor. For each $S\in
\Ext(\sigma)$ we consider the type $\III_1$ factor
$M_{S,1}=M_S\otimes R_\infty$. Since $M_S\otimes R_\infty=M_S\otimes
(R_{0,1}\rtimes_{\sigma^\phi}\R)= (M_S\otimes
R_{0,1})\rtimes_{Id\otimes\sigma^\phi}\R$, the uniqueness of the
continuous decomposition \cite[Theorems XII.1.1 \&
XII.1.7]{takesaki2} implies that for $S_1, S_2 \in \Ext(\sigma)$,
$M_{S_1,1}$ is isomorphic to $M_{S_2,1}$ if and only if their cores
$M_{S_1}\otimes R_{0,1}$ and $M_{S_2}\otimes R_{0,1}$ are
isomorphic. We proceed as in the previous case to conclude that
$\sim_{oe}$ is Borel reducible to isomorphism of $\III_1$ factors.

\bigskip
\noindent {\it The type $\III_0$  case:} Once one fixes an injective
type $\III_0$ factor, one can either use its discrete or continuous
decomposition and repeat the same construction as in the cases
$0<\lambda\leq 1$, to exhibit a family of $\III_0$ factors that are
not classifiable by countable structures. Rather than adapting here
the previous proofs to take into consideration the fact that the
core of a  $\III_0$ factor is not a factor, in the next section we
show instead that isomorphism relation of injective factors is not
classifiable by countable structures.\qed
\end{pf}

\section{The injective type $\III_0$ case}

Foreman and Weiss showed in \cite{foremanweiss} that the isomorphism
relation of ergodic measure preserving transformations on a standard
probability space $(Y,\nu)$ is generically turbulent.  The goal here
is to combine this result with Kreiger's theorem on the
classification of Krieger factors to exhibit a family of injective
type $\III_0$ factors that are not classifiable by countable
structures. In what follows we assume familiarity with the proof of
Krieger's Theorem as presented by Takesaki in
\cite[XVIII.2]{takesaki3} and with the notation used there. Recall
that in this context we say that two measure class preserving
transformations $T_1$ and $T_2$ on a standard measure space are
conjugate if there exists a {\it measure class preserving}
automorphism $\theta$ such that $\theta T_1=T_2 \theta$.

If $Q\in \Aut(Y,\nu)$ is ergodic, we construct the ergodic flow
$F^Q_t$ under the constant ceiling function $r(y)=1$. More
specifically, for each $t\in\R$,
$F^Q_t\in\Aut(Y\times[0,1],\nu\times m)$ is defined as
$F^Q_t(y,s)=(Q^n(y),\alpha)$, where $t+s+\alpha$, $n\in\Z$,
$0\leq\alpha<1$ (see \cite[XII.3]{takesaki2}). To be consistent with
the notation used by Takesaki, we denote with $(\Omega,\mu)$ the
standard probability space $(Y\times[0,1],\nu\times m)$. Since $Q$
is measure preserving, so is the flow $F^Q_t$. Moreover, it can be
shown that two flows $F^Q_t$ and $F^{\widetilde Q}_t$ constructed in
this manner are conjugate if and only if the base transformations
$Q$ and $\widetilde Q$ are conjugate by a {\it measure class
preserving} automorphism, \cite[p.31]{krengel}. Since $Q$ and
$\widetilde Q$ are ergodic and measure preserving, it follows that
they are in fact conjugate by a {\it measure preserving}
automorphism.

Since the map $\R\times\Aut(Y,\nu)\to\Aut( \Omega,\mu):(t,Q)\mapsto
F^Q_t$ is Borel, the map
\begin{align*}
\Aut(Y,\nu)&\to\Aut(X\times \Omega\times \R, P\times\mu \times e^{-s} ds):\\
Q&\mapsto S^Q(x,\omega,r)=(Tx,F^Q_{a(x)}(\omega),r+b(x))
\end{align*}
of \cite[Theorem XVIII.2.5]{takesaki3} is Borel. Moreover,  $S^Q$ is
an ergodic transformation whose associated modular flow is $F^Q_t$.
(Observe that since $F^Q_t$ is m.p., $p(t,\omega)$ defined in
\cite[p.g. 319 (11)]{takesaki3} is equal to $1$.)

By Krieger's theorem (\cite[Theorem XVIII.2.1]{takesaki3}), the type
$\III_0$ factors $L^\infty(X\times \Omega\times \R, P\times\mu
\times e^{-s} ds)\rtimes_{S^Q}\Z$ and $L^\infty(X\times \Omega\times
\R, P\times\mu \times e^{-s} ds)\rtimes_{S^{\widetilde Q}}\Z$ are
isomorphic if and only if the associated modular flows $F^Q_t$ and
$F^{\widetilde Q}_t$ are conjugate. Thus the map $ Q\mapsto
L^\infty(X\times \Omega\times \R, P\times\mu \times e^{-s} ds)
\rtimes_{S^{Q}}\Z$ is a Borel reduction from isomorphism of ergodic
measure preserving transformations to isomorphism of injective
$\III_0$ factors. We have shown:

{\thm The isomorphism relation for injective type $\III_0$ factors
is not classifiable by countable structures.}

\section{Isomorphism is complete analytic}

An equivalence relation $E$ is called {\it Borel complete} for
countable structures if for every language $\mathcal L$ and any
Borel $\mathscr C\subseteq\Mod(\mathcal L)$ invariant under
$\simeq^{\Mod(\mathcal L)}$ we have $\simeq^{\mathscr C}\leq_B E$.
In particular, if $E$ is Borel complete for countable structures
then we can Borel reduce the isomorphism relation of countable
graphs, groups, fields, linear orders, etc., to $E$. In this section
we prove:

{\thm The isomorphism relation $\simeq^{\mathscr F_{\II_1}}$ is
Borel complete for countable structures.\label{thmcompl1}}

In particular, it follows from this and Theorem \ref{basicthm} that
$\simeq^{\gf}<_B\simeq^{\mathscr F_{\II_1}}$ and
$\simeq^{\gp}<_B\simeq^{\mathscr F_{\II_1}}$. Since it is known that
$\simeq^{\gf}$ and $\simeq^{\gp}$ are complete analytic as subsets
of $\gf\times\gf$ and $\gp\times\gp$, respectively, we obtain from
Theorem \ref{thmcompl1} that $\simeq^{\mathscr F_{\II_1}}$ is a
complete analytic subset of $\mathscr F_{\II_1}\times\mathscr
F_{\II_1}$, see Corollary \ref{corana} below.

Let $\wicc$ denote the subset of $\gp$ consisting of ICC groups with
the relative property (T) over an infinite normal subgroup. It may
be verified that $\wicc$ is a Borel subset of $\gp$. To prove
Theorem \ref{thmcompl1} we will use a deep result of Sorin Popa
\cite[Theorem 7.1]{popa3} which implies that for groups
$G_1,G_2\in\wicc$ it holds that if the $\II_1$ factors $M_{G_1}$ and
$M_{G_2}$ that arise from Bernoulli shifts of $G_1$ and $G_2$ via
the group measure space construction are isomorphic then the groups
$G_1$ and $G_2$ are isomorphic.

{\lem The equivalence relation $\simeq^{\wicc}$ is Borel complete
for countable structures.}

\begin{pf}
By a graph we mean a symmetric, irreflexive relation. Graphs are not
assumed to be connected unless otherwise stated. All graphs here are
countable.

We will modify Mekler's construction in \cite{mekler}. Recall that
Mekler introduces a notion of ``nice'' graph \cite[Defn.
1.1]{mekler}, and shows in \cite[\S 2]{mekler} that for each nice
graph $\Gamma$ there is an associated nil-2 exponent $p$ group
$G(\Gamma)$, which we will call the Mekler group of $\Gamma$. As
noted in \cite{friedstan}, the association $\Gamma\mapsto G(\Gamma)$
is a Borel reduction of isomorphism of nice graphs to isomorphism of
nil-2 exponent $p$ groups. Mekler also shows that the isomorphism
relation of {\it connected} nice graphs is Borel complete for
countable structures, \cite[\S 1]{mekler}.

Let $\Gamma$ be a connected nice graph. Define a graph
$\Gamma_{\F_2}$ on $\F_2\times\Gamma$ by
$$
(a_1,v_1)\Gamma_{\F_2}(a_2,v_2)\iff a_1=a_2\wedge v_1\Gamma v_2.
$$
Then $\Gamma_{\F_2}$ is nice (but not connected). The shift action
of $\sigma:\F_2\curvearrowright\Gamma_{\F_2}$ defined by
$$
\sigma(a)(b,v)=(a^{-1}b,v)
$$
is clearly an action by graph automorphisms of $\F_2$ on
$\Gamma_{\F_2}$. Recall that for a fixed prime $p$, the associated
Mekler group $G(\Gamma_{\F_2})$ is defined as
$$
\big (\underset{(a,v)\in \Gamma_{\F_2}}{\mathbf 2} \Z/p\Z\big )/N
$$
where $N=\langle [(a_1,v_1),(a_2,v_2)]:
(a_1,n_1)\Gamma_{\F_2}(a_2,n_2)\rangle$, and $\mathbf 2$ denotes the
free product in the category of nil-2 exponent $p$ groups. The
action $\sigma$ corresponds to an action of
$\sigma_0:\F_2\curvearrowright\underset{(a,v)\in
\Gamma_{\F_2}}{\mathbf 2} \Z/p\Z$ by group automorphisms in the
obvious way. Since $\sigma$ acts by graph automorphisms, $N$ is
$\sigma_0$-invariant. Hence $\sigma_0$ factors to an action
$\hat\sigma$ of $\F_2$ on $G(\Gamma_{\F_2})$ by group automorphisms,
and we may form the semi-direct product
$G(\Gamma_{\F_2})\rtimes_{\hat\sigma} \F_2$.

{{\bf Claim.} $G(\Gamma_{\F_2})\rtimes_{\hat\sigma} \F_2$ is ICC}.

\medskip

\noindent {\it Proof.} Since $\F_2$ is ICC, it is enough to consider
elements of the form $(g,e)\in G(\Gamma_{\F_2})\rtimes_{\hat\sigma}
\F_2$, $g\neq e$. Since $(e,a)(g,e)(e,a^{-1})=(\hat\sigma(a)(g),e)$
it is enough to show that
$$
[g]_{\hat\sigma}=\{\hat \sigma(a)(g):a\in\F_2\}
$$
is infinite. Suppose for a contradiction that this set is finite.
Let
$$
\mathbf 2 \Gamma_{\F_2}= \underset{(a,v)\in \Gamma_{\F_2}}{\mathbf
2} \Z/p\Z
$$
and identify this group with the free nil-2 exponent $p$ group
generated by the vertices of the graph $\Gamma_{\F_2}$. Let
$\varphi: \mathbf 2\ \Gamma_{\F_2}\to G(\Gamma_{\F_2})$ be the
quotient map with $\ker(\varphi)=N$.

If $[g]_{\hat\sigma}$ is finite, then we may find finitely many
nodes $(a_1,v_1),\ldots (a_k,v_k)\in\Gamma_{\F_2}$ such that each
element of $[g]_{\hat\sigma}$ can be written as
$$
\varphi((a_{i_1},v_{i_1})\cdots(a_{i_l},v_{i_l}))
$$
for an appropriate choice of $1\leq i_1,\ldots, i_l\leq k$. Let
$\Delta=\{a_{1},\ldots, a_{k}\}$ and let $b\in\F_2$ be such that
$b^{-1}\Delta\cap\Delta=\emptyset$. Since $\mathbf 2 \Gamma_{\F_2}$
is freely generated, we have a natural homomorphism
$$
p_{\Delta}:\mathbf 2 \Gamma_{\F_2}\to \mathbf 2 \Delta\times\Gamma,
$$
under which
$$
p_{\Delta}((a_{i_1},v_{i_1})\cdots(a_{i_l},v_{i_l}))=(a_{i_1},v_{i_1})\cdots(a_{i_l},v_{i_l}).
$$
Then since
$$
\sigma_0(b)((a_{i_1},v_{i_1})\cdots(a_{i_l},v_{i_l}))=(b^{-1}a_{i_1},v_{i_1})\cdots(b^{-1}a_{i_l},v_{i_l})
$$
we have
$p_{\Delta}(\sigma_0(b)((a_{i_1},v_{i_1})\cdots(a_{i_l},v_{i_l})))=e$.
It follows that if $\hat\sigma(b)(g)$ is represented by a reduced
word in $(a_1,v_1),\ldots, (a_k,v_k)$ it must be the empty word.
Hence $\hat\sigma(b)(g)=e$ and so $g=e$, a contradiction. \hfill
[Claim]\qed

Define
$$
G_{\Gamma}=\SL(3,\Z)\times G(\Gamma_{\F_2})\rtimes_{\hat\sigma}\F_2.
$$
Since $\SL(3,\Z)$ has property (T) of Kazhdan and is ICC, then
$G_{\Gamma}\in \wicc$. We claim that the map $\Gamma\mapsto
G_\Gamma$ is a Borel reduction of isomorphism of nice connected
graphs to $\simeq^{\wicc}$. It is routine to verify that
$\Gamma\mapsto G_\Gamma$ is Borel. It is also clear that if
$\Gamma\simeq\Gamma'$ then $G_\Gamma\simeq G_{\Gamma'}$. So suppose
that $G_{\Gamma}\simeq G_{\Gamma'}$. Let
$$
H_0=\langle g\in G_{\Gamma}: \,(\exists\chi\in \Char(G_\Gamma))\
\chi(g)\neq 1\rangle.
$$
It is clear that $\F_2\subseteq H_0\subseteq
G(\Gamma_{\F_2})\rtimes_{\hat\sigma}\F_2$, since $\SL(3,\Z)$ has no
non-trivial characters. Since $G(\Gamma_{\F_2})$ is an exponent $p$
group we have
$$
G(\Gamma_{\F_2})=\{g\in G_\Gamma: g=e\vee (g^p=e\wedge(\exists h\in
H_0) hgh^{-1}\neq g)\}
$$
from which it follows that if $G_{\Gamma}\simeq G_{\Gamma'}$ then
$$
G(\Gamma_{\F_2})\simeq G(\Gamma'_{\F_2}).
$$
By \cite[Lemma 2.2]{mekler} it follows that
$\Gamma_{\F_2}\simeq\Gamma'_{\F_2}$. Since the connected components
of $\Gamma_{\F_2}$ and $\Gamma'_{\F_2}$ are all isomorphic to
$\Gamma$ and $\Gamma'$, respectively, it follows that
$\Gamma\simeq\Gamma'$.

As noted above, the isomorphism relation of connected nice graphs is
Borel complete for countable structures, and so it follows that
$\simeq^{\wicc}$ is also Borel complete for countable structures.
\qed \end{pf}

{\thm It holds that $\simeq^{\wicc}<_B\simeq^{\mathscr F_{\II_1}}$.
Hence $\simeq^{\mathscr F_{\II_1}}$ is Borel complete for countable
structures, but not classifiable by countable
structures.\label{thmcompl2}}

\medskip

It is clear that Theorem \ref{thmcompl2} implies Theorem
\ref{thmcompl1}.

\medskip

\noindent {\it Notation}: If $G\in\gp$, we will write $\cdot_G$ for
composition in $G$, $^{-1_G}$ for the inverse operation in $G$, and
$e_G$ for the identity in $G$.

\begin{pf}[Proof of Theorem \ref{thmcompl2}]
Let $X=[0,1]^{\N}$, equipped with the product measure $\mu$. For
each $G\in\gp$ let $\sigma_{G}:G\curvearrowright X$ be the Bernoulli
shift action,
$$
\sigma_G(g)(x)(g_0)=x(g^{-1_G}\cdot_G g_0).
$$
For each $G$, let $M_G=L^\infty(X)\rtimes_{\sigma_G} G$.

{{\bf Claim.} For each $g\in\N$, the map $G\mapsto \sigma_G(g)$ is a
continuous map $\gp\to\Aut(X,\mu)$.}

\medskip

\noindent{\it Proof of Claim.} Fix intervals $I_1,\ldots,
I_k\subseteq [0,1]$ and $n_1,\ldots,n_k\in\N$. Consider the cylinder
set
$$
B=\{x\in X: (\forall i\leq k) x(n_i)\in I_i\},
$$
and fix $G_0\in\gp$. Define
$$
N(G_0)=\{G\in\gp: (\forall i\leq k) g\cdot_G n_i=g\cdot_{G_0} n_i\}.
$$
Clearly $N(G_0)$ is a clopen subset of $\gp\subseteq
\N^{\N\times\N}\times\N$ in the subspace topology. If $G\in N(G_0)$
then
\begin{align*}
x\in\sigma_G(g)(B)&\iff \sigma_G(g^{-1_G})(x)\in B\\
&\iff (\forall i\leq k) \sigma_G(g^{-1_G})(x)(n_i)\in I_i\\
&\iff (\forall i\leq k) x(g\cdot_G n_i)\in I_i\\
&\iff (\forall i\leq k) x(g\cdot_{G_0} n_i)\in I_i\\
&\iff (\forall i\leq k)
\sigma_{G_0}(g^{-1_{G_0}})(x)(n_i)\in I\\
&\iff x\in\sigma_{G_0}(g)(B)
\end{align*}
Thus for $G\in N(G_0)$ we have $\sigma_G(g)(B)=\sigma_{G_0}(g)(B)$.
Since any measurable $C\subseteq X$ may be approximated by a finite
union of cylinder sets, it follows from the above that
$G\mapsto\sigma_G(g)$ is continuous. \hfill [Claim]\qed

\medskip

Using the claim, it is clear that for each $g\in\N$ the map
$G\mapsto U^{\sigma_G}_{g}$, where $U^{\sigma_G}_{g}$ is the unitary
operator on $L^2(X\times \N)$ defined by
$$
U^{\sigma_G}_{g}(f)(x,g_0)=f(\sigma_G(g)(x),g^{-1_G}\cdot_G g_0),
$$
is continuous. Exactly as in Lemma \ref{basiclemma} we then have
$$
G\mapsto M_G=L^\infty(X)\rtimes_{\sigma_G} G
$$
is Borel, since $M_G$ is generated by $U^{\sigma_G}_{g}$, $g\in G$,
and $L_\psi$, $\psi\in L^\infty(X)$. We claim that the map $\wicc\to
\mathscr F_{\II_1}:G\mapsto M_G$ is a Borel reduction of
$\simeq^{\wicc}$ to $\simeq^{\mathscr F_{\II_1}}$. Clearly, if
$G_0,G_1\in\wicc$ are isomorphic then $\sigma_{G_0}$ and
$\sigma_{G_1}$ are conjugate, and so $M_{G_0}$ and $M_{G_1}$ are
isomorphic. On the other hand, if $M_{G_0}\simeq M_{G_1}$ then we
may apply Theorem 7.1 in \cite{popa3} to conclude that $G_0\simeq
G_1$. This together with Theorem \ref{basicthm} proves that
$\simeq^{\wicc}<_B\simeq^{\mathscr F_{\II_1}}$.\qed
\end{pf}

{\it Remark.}  It should be noted that Corollary 7.2 in \cite{popa3}
subsumes that at least one of the groups in question is ICC. Thus it
does {\it not} suffice for Theorem \ref{thmcompl2} to look at groups
of the form $\SL(3,\Z)\times G(\Gamma)$ for some nice graph
$\Gamma$.

{\cor The isomorphism relation $\simeq^{\mathscr F_{\II_1}}$ is a
complete analytic subset of $\mathscr F_{\II_1}\times\mathscr
F_{\II_1}$. In particular, it is not Borel.\label{corana}}

\begin{pf}
Let $\simeq^{\mathscr F}_u$ denote unitary equivalence in $\mathscr
F$. By \cite[Lemma 2.1]{effros2} the action of $\mathcal U(\mathcal
H)$ on $\mathscr F$ is Borel, thus $\simeq^{\mathscr F}_u$ is
analytic. As noted on p. 436 in \cite{effros2}, the map
$$
M\mapsto M\otimes \C\in \mathscr F(\mathcal H\otimes\mathcal H),
$$
constitutes a Borel reduction of $\simeq^{\mathscr F(\mathcal H)}$
to $\simeq_u^{\mathscr F(\mathcal H\otimes \mathcal H)}$, since any
algebraic isomorphism between two separable von Neumann algebras
with properly infinite commutants is unitarily implemented. It
follows that $\simeq^{\mathscr F}$ is analytic. To see that it is
complete analytic, note that by the previous Theorem the isomorphism
relation of countable Abelian $p$-groups is Borel reducible to
$\simeq^{\mathscr F_{\II_1}}$. Since this equivalence relation is
complete analytic by \cite[Theorem 6]{friedstan}, it follows that
$\simeq^{\mathscr F_{\II_1}}$ is complete analytic.\qed
\end{pf}

{\it Remark.}  It is possible to derive Corollary \ref{corana} as a
consequence of \cite[Corollary 0.7]{iopepo}, by showing that the
construction given there of a factor $M$ with $\Out(M)\simeq K$ for
a prescribed compact Abelian group $K$ is Borel in the codes.
However, it is not known if isomorphism of countable Abelian groups
is Borel complete for countable structures, so this approach does
not give a different proof of Theorem \ref{thmcompl2}.

\paragraph{} While the results of this paper indicate that the isomorphism
relation for factors is, in general, very complicated, we point out
that it is not as complicated (from the point of view of $\leq_B$)
as some other analytic equivalence relations. An equivalence
relation $E$ on a standard Borel space is {\it below a group action}
if there is a Borel action of a Polish group $G$ on a Polish space
$X$ such that $E\leq_B E_G$, where $E_G$ denotes the orbit
equivalence relation in $X$. Then we have:

{\thm {\rm (}i{\rm )} Isomorphism of separable von Neumann factors
is below a group action.

{\rm (}ii{\rm )} If $H$ is a countable group then orbit equivalence
of probability measure preserving a.e. free ergodic actions of $H$
is below a group action.}

\begin{pf}
We have already seen in the proof of Corollary \ref{corana} that
$$
M\mapsto M\otimes \C\in \mathscr F(\mathcal H\otimes\mathcal H),
$$
is a Borel reduction of algebraic isomorphism to unitary
equivalence, hence ({\it i}) is clear.

For ({\it ii}), let $(X,\mu)$ be a standard Borel probability space
and let $\mathcal H=L^2(H\times X)$. Let
$$
\mathcal U_{L^\infty(X)\otimes \C}(\mathcal H\otimes\mathcal
H)=\{u\in \mathcal U(\mathcal H\otimes\mathcal H):
uL^{\infty}(X)\otimes\C u^*=L^{\infty}(X)\otimes\C\}
$$
be the stabilizer of $L^\infty(X)\otimes \C$, which is a closed
subgroup of $\mathcal U(\mathcal H\otimes\mathcal H)$. Let $E$
denote the orbit equivalence relation induced by the conjugation
action of $\mathcal U_{L^\infty(X)\otimes \C}(\mathcal
H\otimes\mathcal H)$ on $\vN(\mathcal H\otimes\mathcal H)$. We claim
that orbit equivalence of a.e. free ergodic actions of $H$ on
$(X,\mu)$ is Borel reducible to $E$.

For $\sigma\in\act(H,X,\mu)$ let
$M_\sigma=L^\infty(X)\rtimes_{\sigma} H$. Then the map
$$
\act(H,X,\mu)\to\vN(\mathcal H\otimes\mathcal H):\sigma\mapsto
M_{\sigma}\otimes\C
$$
is Borel, and if $\sigma_0,\sigma_1\in\act(H,X,\mu)$ are orbit
equivalent a.e. free ergodic actions then the inclusions
$L^\infty(X)\subset M_{\sigma_0}$ and $L^\infty(X)\subset
M_{\sigma_1}$ are isomorphic. It follows that there is a unitary
$u\in\mathcal U_{L^\infty(X)\otimes \C}$ such that
$$
uM_{\sigma_0}\otimes\C u^*=M_{\sigma_1}\otimes\C.
$$
If conversely there is a unitary $u\in\mathcal U_{L^\infty(X)\otimes
\C}$ such that
$$
uM_{\sigma_0}\otimes\C u^*=M_{\sigma_1}\otimes\C
$$
then we have the isomorphic inclusions
\begin{align*}
L^\infty(X)\subset M_{\sigma_0}&\simeq L^\infty(X)\otimes\C\subset M_{\sigma_0}\otimes\C\\
&\simeq L^\infty(X)\otimes\C\subset M_{\sigma_1}\otimes\C\simeq
L^\infty(X)\subset M_{\sigma_1},
\end{align*}
and so by Feldman and Moore's Theorem, $\sigma_0$ is orbit
equivalent to $\sigma_1$. \qed \end{pf}

{\it Remark.}  By \cite[Theorem 4.2]{kechlouv}, it follows from the
previous Theorem that the equivalence relation $E_1$ on $\R^\N$,
defined by
$$
xE_1y\iff (\exists N)(\forall n\geq N) x(n)=y(n)
$$
is not Borel reducible to $\simeq^{\mathscr F}$. In particular,
$\simeq^{\mathscr F}$ is not universal for analytic relations, i.e.,
not every analytic equivalence relation is Borel reducible to
$\simeq^{\mathscr F}$. Part ({\it ii}) of the previous Theorem
answers a question raised by Hjorth.

\bigskip

\begin{small}
{\sc\noindent  Instituto de Ciencias\\
Universidad Nacional de General Sarmiento\\
J. M. Gutierrez 1150\\
(1613) Los Polvorines, Argentina

\medskip

Instituto Argentino  de Matem\'aticas-CONICET\\
Saavedra 15, Piso 3
(1083), Buenos Aires, Argentina}\\
{\it E-mail}: {\tt rsasyk@ungs.edu.ar}
\end{small}

\bigskip
\begin{small}
{\sc\noindent Kurt G\"odel Research Center\\
University of Vienna\\
W\"ahringer Strasse 25\\
1090 Vienna, Austria\\
 {\it E-mail}: {\tt asger@logic.univie.ac.at}}
\end{small}
\end{document}